\documentclass[10pt,twoside]{article}
\usepackage{amssymb,amscd,amsmath}
\usepackage{graphicx}
\usepackage{Latex-document}

\newtheorem{theorem}{Theorem}
\newtheorem{definition}[theorem]{Definition}

\newtheorem{thmnn}{Theorem} 
\newtheorem{cornn}{Corollary} 

\def\C{\mathbb C}
\def\K{\mathbb K}
\def\N{\mathbb N}
\def\Q{\mathbb Q}
\def\R{\mathbb R}
\def\Z{\mathbb Z}

\def\cA{\mathcal A}
\def\cB{\mathcal B}
\def\cG{\mathcal G}

\def\cN{\mathcal N}

\newcommand{\hra}{\hookrightarrow}

\title{\bf Knots, von Neumann Signatures, \vskip -2mm and Grope
Cobordism\thanks{Partially supported by an NSF-grant and the
Max-Planck Gesellschaft.}  \vskip 6mm}

\author{Peter Teichner\vspace*{-0.5cm}\thanks{University of California in
San Diego, 9500 Gilman Drive, La Jolla, CA 92093-0112, USA.
E-mail: teichner@math.ucsd.edu}}

\date{\vspace{-8mm}}

\begin{document}

\maketitle

\thispagestyle{first} \setcounter{page}{437} \setcounter{figure}{0}

\begin{abstract}

\vskip 3mm

We explain new developments in classical knot theory in 3 and
4~dimensions, i.e.\ we study knots in 3-space, up to isotopy as well as up
to concordance. In dimension~3 we give a geometric interpretation of the
Kontsevich integral (joint with Jim Conant), and in dimension~4 we
introduce new concordance invariants using von Neumann signatures (joint
with Tim Cochran and Kent Orr). The common geometric feature of our
results is the notion of a grope cobordism.
\vskip 4.5mm



\end{abstract}
\vskip 12mm

\section{Introduction}
\setzero\vskip-5mm \hspace{5mm}

A lot of fascinating mathematics has been created when successful tools are
transferred from one research area to another. We shall describe two
instances of such transfers, both into knot theory. The first transfer
realizes commutator calculus of group theory by {\em embedded versions}
in 3- and 4-space, and produces many interesting geometric
equivalence relations on knots, called {\em grope cobordism} in 3-space
and {\em grope concordance} in 4-space. It turns out that in 3-space
these new equivalence relations give a geometric interpretation
(Theorem~\ref{thm:3d}) of Vassiliev's finite type invariants
\cite{V} and that the Kontsevich integral \cite{K} calculates the new
theory over $\Q$ (Theorem~\ref{thm:kontsevich}).

In 4-space the new equivalence relations
factor naturally through knot concordance, and in fact they organize all
known concordance invariants in a wonderful manner (Theorem~\ref{thm:4d}).
They also point the way to new concordance invariants
(Theorem~\ref{thm:cot2}) and these are constructed using a
second transfer, from the spectral theory of self-adjoint operators and
von Neumann's continuous dimension \cite{MvN}.

\subsection{A geometric interpretation of group commutators}

\hspace{5mm} To explain the first transfer into knot theory,
recall that every knot bounds a Seifert surface (embedded in
3-space), but only the trivial knot bounds an embedded disk. Thus
all of knot theory is created by the difference between a surface
and a disk. The new idea is to filter this difference by
introducing a concept into knot theory which is the analogue of
iterated commutators in group theory. Commutators arise because a
continuous map $\phi:S^1\to X$ extends to a map of a surface if
and only if $\phi$ represents a commutator in the fundamental
group $\pi_1X$. Iterated commutators can similarly be expressed by
gluing together several surfaces. Namely, there are certain finite
2-complexes (built out of iterated surface stages) called {\em
gropes} by Cannon \cite{C}, with the following defining property:
$\phi:S^1\to X$ represents an element in the $k$-th term of the
{\em lower central series} of $\pi_1X$ if and only if it extends
to a continuous map of a {\em grope of class~$k$}. Similarly,
there are {\em symmetric gropes} which geometrically implement the
derived series of $\pi_1X$, see Figures~\ref{fig:gropes} and
\ref{fig:symmetric gropes}

Gropes, therefore, are not quite manifolds but the singularities that arise
are of a very simple type, so that these 2-complexes are in some sense the
next easiest thing after surfaces. Two sentences on the history of the use
of gropes in mathematics are in place, compare \cite[Sec.2.11]{FQ}. Their
inventor Stan'ko worked in high-dimensional topology, and so did Edwards
and Cannon who developed gropes further. Bob Edwards suggested their
relevance for topological 4-manifolds, where they were used extensively,
see
\cite{FQ} or \cite{FT}. It is this application that seems to have created a
certain ''Angst`` of studying gropes, so we should point out that the only
really difficult part in that application is the use of {\em infinite
constructions}, i.e.\ when the class of the grope goes to infinity.

One purpose of this note is to
convince the reader that (finite) gropes are a very simple, and extremely
powerful tool in low-dimensional topology. The point is that once one can
describe iterated commutators in $\pi_1X$ by maps of gropes, one might as
well study {\em embedded} gropes in 3-space (respectively 4-space) in order
to organize knots up to isotopy (respectively up to concordance).
In Section~\ref{sec:3d} we shall explain joint work
with Jim Conant on how gropes embedded in 3-space
lead to a geometric interpretation of Vassiliev's knot invariants \cite{V}
and of the Kontsevich integral \cite{K}.

\subsection{von Neumann signatures and knot concordance}

\hspace{5mm} In Section~\ref{sec:4d} we study symmetric gropes
embedded in 4-space, and explain how they lead to a geometric
framework for all known knot concordance invariants and beyond.
More precisely, we explain our joint work with Tim Cochran and
Kent Orr \cite{COT1}, where we define new concordance invariants
by inductively constructing representations into certain  solvable
groups $G$, and associating a hermitian form over the group ring
$\Z G$ to the knot $K$, which is derived from the intersection
form of a certain $4$-manifold with fundamental group $G$ and
whose boundary is obtained by $0$-surgery on $K$.  This
intersection form represents an element in the Cappell-Shaneson
$\Gamma$-group \cite{CS} of $\Z G$  and we detect it via the
second transfer from a different area of mathematics: The standard
way to detect elements in Witt groups like the $\Gamma$-group
above is to construct unitary representations of $G$, and then
consider the corresponding (twisted) signature of the resulting
hermitian form over $\C$. It turns out that the solvable groups
$G$ we construct do not have any interesting finite dimensional
representations, basically because they are ``too big'' (e.g.\ not
finitely generated), a property that is intrinsic to the groups
$G$ in question because they are ``universal solvable'' in the
sense that many $4$-manifold groups (with the right boundary) must
map to $G$, extending the given map of the knot group.

However, every group $G$ has a fundamental
unitary representation given by $\ell^2G$, the Hilbert space of square
summable sequences of group elements with complex coefficients. The
resulting (weak) completion of $\C G$ is the {\em group von Neumann
algebra} $\cN G$. It is of type~$II_1$ because the map $\sum a_g g\mapsto
a_1$ extends from $\C G$ to give a finite faithful trace on $\cN G$.

The punchline is that hermitian forms over the completion $\cN G$ are much
easier to understand than over $\C G$ because they are diagonalizable
(by functional calculus of self-adjoint operators). Here one really uses
the von Neumann algebra, rather than the $C^*$-algebra completion of $\C G$
because the functional calculus must be applied to the characteristic
functions of the positive (respectively negative) real numbers, which are
bounded but {\em not} continuous.

The subspace on which the hermitian form is positive
(respectively negative) definite has a {\em continuous dimension}, which
is the positive real number given by the trace of the projection onto that
subspace. As a consequence, one can associate to every hermitian form over
$\cN G$ a real valued invariant, the {\em von Neumann signature}.
In \cite{COT1} we use this invariant to construct our new knot
concordance invariants, and a survey of this work can be found in
Section~\ref{sec:4d} It is not only related to embedded gropes in 4-space
but also to the existence of towers of {\em Whitney disks} in 4-space.
Unfortunately, we won't be able to explain this aspect of the theory, but
see \cite[Thm.8.12]{COT1}.

\subsection{Noncommutative Alexander modules}

\hspace{5mm} In Section~\ref{sec:4d} we shall hint at how the
interesting representations to our solvable groups are obtained.
But it is well worth pointing out that the methods developed for
studying knot concordance have much simpler counterparts in
3-space, i.e.\ if one is only interested in isotopy invariants.

A typical list of knot invariants that might find its way into a text book
or survey talk on {\em classical} knot theory, would  contain the
Alexander polynomial, (twisted) signatures, (twisted) Arf invariants, and
maybe knot determinants. It turns out that all of these invariants can be
computed from the homology of the infinite cyclic covering of the knot
complement, and are in this sense ``commutative'' invariants.

Instead of the maximal abelian quotient one can use other solvable quotient
groups of the knot group to obtain ``noncommutative'' knot invariants.
The canonical candidates are the quotient groups $G_n$ of the derived
series $G^{(n)}$ of the knot group (compare Section~\ref{sec:4d} for the
definition). One can thus define the higher order {\em Alexander modules}
of a knot $K$ to be the $\Z G_{n+1}$-modules
$$
\cA_{n}(K):=H_1(S^3 \smallsetminus K;\Z G_{n+1}).
$$
The indexing is chosen so that $\cA_0$ is the classical
Alexander module. For $n\geq 1$ these modules are best studied by
introducing further algebraic tools as follows. By a result of Strebel the
groups $G_n$ are torsionfree. Therefore, the group ring
$\Z G_n$ satisfies the Ore condition and has a well defined (skew) quotient
field. This field is in fact the quotient field of a (skew) polynomial
ring $\K_n[t^{\pm 1}]$, with $\K_n$ the quotient field of
$\Z[G^{(1)}/G^{(n)}]$ and $G_1= \langle t\rangle
\cong\Z$. Thus one is exactly in the context of \cite[Sec.2]{COT1} and one
can define explicit noncommutative isotopy invariants of knots. For
example, let $d_n(K)$ be the dimension (over the field $\K_{n+1}$) of the
{\em rational} Alexander module
$$
\cA_n(K) \otimes_{\Z G_{n+1}} \K_{n+1}[t^{\pm 1}].
$$
It is shown in \cite[Prop.2.11]{COT1} that these dimensions are finite with
the degree of the usual Alexander polynomial being $d_0(K)$. Moreover,
Cochran \cite{Co} has proven the following non-triviality result for these
dimensions.

\begin{thmnn}
If $K$ is a nontrivial knot then for $n\geq 1$ one has
$$
d_0(K)\leq d_1(K)+1\leq d_2(K)+1\leq \dots \leq d_n(K)+1 \leq
2\cdot\text{genus of }K.
$$
Moreover, there are examples where these numbers are strictly
increasing up to any given $n$.
\end{thmnn}

\begin{cornn}
If one of the inequalities in the above theorem is strict then $K$
is not fibered. Furthermore, $0$-surgery on $K$ cross the circle is not a
symplectic 4-manifold.
\end{cornn}

The first statement is clear: For fibered knots the
degree of the Alexander polynomial $d_0(K)$ equals twice the genus of the
knot $K$. The second statement follows from a result of Kronheimer
\cite{Kr} who showed that this equality also holds if the above
$4$-manifold is symplectic.

Recently, Harvey \cite{Ha} has studied similar invariants for arbitrary
$3$-manifolds and has proven generalizations of the above results:
There are lower bounds for the Thurston norm
of a homology class, analogous to $d_i(K)$, that are better than
McMullen's lower bound, which is the analogy of $d_0(K)$. As a
consequence, she gets  new algebraic obstructions to a $4$-manifold of the
form $M^3\times S^1$ admitting a symplectic structure.

Just like in the classical case $n=0$, there is more structure on the
rational Alexander modules. By \cite[Thm.2.13]{COT1} there are higher order
{\em Blanchfield forms} which are hermitian and non-singular in an
appropriate sense, compare \cite[Prop.12.2]{Co}.
It would be very interesting to know whether the $n$-th order
Blanchfield form determines the von Neumann $\eta$-invariant associated to
the $G_{n+1}$-cover. So far, these $\eta$-invariants are very mysterious
real numbers canonically associated to a knot.

Only in the bottom case $n=0$ do we understand this $\eta$-invariant well:
The $L^2$-index theorem implies that the von Neumann $\eta$-invariant
corresponding to the $\Z$-cover is the von Neumann signature of a certain
$4$-manifold with fundamental group $\Z$. Moreover, this signature is the
integral, over the circle, of all (Levine-Tristram) twisted signatures of
the knot \cite[Prop.5.1]{COT2} (and is thus a concordance invariant). For
$n\geq 1$ there is in general no such $4$-manifold available and the
corresponding $\eta$-invariants are not concordance invariants.

\section{Grope cobordism in 3-space} \label{sec:3d}
\setzero\vskip-5mm \hspace{5mm}

We first give a more precise treatment of the first transfer from
group theory to knot theory hinted at in the introduction. Recall that the
fundamental group consists of continuous maps of the circle $S^1$ into some
target space $X$, modulo homotopy (i.e.\
$1$-parameter families of continuous maps). Quite analogously, classical
knot theory studies smooth {\em embeddings} of a circle into
$S^3$, modulo isotopy (i.e.\ $1$-parameter families of embeddings).
To explain the transfer, we recall that a continuous map $\phi:S^1\to X$
represents the trivial element in the fundamental group $\pi_1X$ if and
only if it extends to a map of the disk, $\tilde\phi:D^2\to X$. Moreover,
$\phi$ represents a commutator in
$\pi_1X$ if and only if it extends to a map of a surface (i.e.\ of
a compact oriented 2-manifold with boundary $S^1$). The first statement has
a straightforward analogy in knot theory: $K:S^1\hra S^3$ is trivial if and
only if it extends to an embedding of the disk into $S^3$. However, every
knot ``is a commutator'' in the sense that it bounds a {\em Seifert
surface}, i.e.\ an embedded surface in $S^3$.

\begin{figure}[ht] \label{fig:gropes}
\begin{center}
\includegraphics[scale=.5]{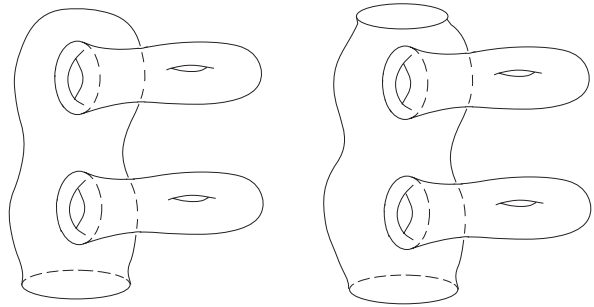}
\caption{Gropes of class~3, with one respectively two boundary
circles}
\end{center}
\end{figure}

Recall from the introduction that gropes are finite 2-complexes defined by
the following property:
$\phi:S^1\to X$ represents an element in the $k$-th term $\pi_1X_k$ of the
lower central series of $\pi_1X$ if and only if it extends to a
continuous map of a grope of class~$k$. Here $G_k$ is defined
inductively for a group $G$ by the iterated commutators
$$
G_2:=[G,G] \text{ and } G_k:=[G,G_{k-1}] \text{ for } k>2.
$$
Accordingly, a grope of class~2 is just a surface, and one can obtain a
grope of class~$k$ by attaching gropes of class~$(k-1)$ to $g$
disjointly embedded curves in the {\em bottom} surface. Here $g$ is the
genus of the bottom surface and the curves are assumed to span one half of
its homology. This gives gropes of class~$k$ with one boundary circle as on
the left of Figure~\ref{fig:gropes} It's not the most general way to get
gropes because of re-bracketing issues, and we refer to
\cite[Sec.2.1]{CT1} for details. The boundary of a grope is by definition
just the boundary of the bottom surface, compare Figure~\ref{fig:gropes}

\begin{definition}  \label{def:3d}
Two (smooth oriented) knots in $S^3$ are {\em grope
cobordant} of class~$k$, if there is an embedded grope of class~$k$ in
$S^3$ (the {\em grope cobordism}) such that its boundary consists exactly
of the given knots.
\end{definition}

An {\em embedding} of a grope is best defined via the obvious $3$-dimensional local model. Since every grope has a
1-dimensional spine, embedded gropes can then be isotoped into the neighborhood of a 1-complex. As a consequence,
embedded gropes abound in 3-space! It is important to point out that if two knots $K_i$ cobounds a grope then
$K_1$ and $K_2$ might very well be linked in a nontrivial way. Thus it is a much stronger condition on $K$ to
assume that it is the boundary of an embedded grope than to say that it cobounds a grope with the unknot. For
example, if $K$ bounds an embedded grope of class~$3$ in $S^3$ then the Alexander polynomial vanishes. Together
with Stavros Garoufalidis, we recently showed \cite{GT} that the 2-loop term of the Kontsevich integral detects
many counterexamples to the converse of this statement.

In joint work with Jim Conant \cite{CT1}, we show that grope cobordism
defines equivalence relations on knots, one for every
class~$k\in\N$. Moreover, Theorem~\ref{thm:3d} below implies that
the resulting quotients are in fact finitely generated abelian groups
(under the connected sum operation). For the smallest values $k=2,3,4$ and
5, these groups are isomorphic to
$$
\{0\}, \Z/2, \Z \text{ and } \Z \times \Z/2
$$
and they are detected by the first two Vassiliev invariants
\cite[Thm.4.2]{CT2}.

The following theorem is formulated in
terms of {\em clasper surgery} which was introduced
independently by Habiro \cite{H} and Goussarov \cite{G}, as a geometric
answer to finite type invariants a l\'a Vassiliev \cite{V}. We cannot
explain the definitions here but see \cite[Thm.1 and 3]{CT1}. We should say
that the notion of a {\em capped} grope is well known in 4~dimensions, see
\cite[Sec.2]{FQ}. In our context, it means that all circles at the ``tips''
of the grope bound disjointly embedded disks in 3-space which are only
allowed to intersect the bottom surface of the grope.

\begin{theorem} \label{thm:3d}
Two knots $K_0$ and $K_1$ are grope cobordant of class~$k$ if and only if
$K_1$ can be obtained from $K_0$ by a finite sequence of clasper surgeries
of {\em grope degree~$k$} (as defined below).

Moreover, two knots are {\em capped} grope cobordant of class~$k$ if
and only if they have the same finite type invariants of Vassiliev
degree~$<k$.
\end{theorem}

As a consequence of this result, the invariants associated to
grope cobordism are highly nontrivial as well as manageable. For example,
we prove the following result in \cite[Thm.1.1]{CT2}:

\begin{theorem} \label{thm:kontsevich}
The (logarithm of the) Kontsevich integral (with values in
$\cB^g_{<k}$), graded by the new {\em grope degree~$k$}, is an obstruction
to finding a grope cobordism of class~$k$ between two knots. Moreover,
this invariant is rationally faithful and surjective.
\end{theorem}

Here $\cB^g_{<k}$ is one of the usual algebras of Feynman diagrams known
from the theory of finite type invariants, but graded by the {\em grope
degree}. More precisely, $\cB^g_{<k}$ is the $\Q$-vector space generated by
connected uni-trivalent graphs of grope degree~$i, 1<i<k$, with at
least one univalent vertex and a cyclic ordering at each trivalent vertex.
The relations are the usual IHX and AS relations. The grope degree is the
Vassiliev degree (i.e.\ half the number of vertices) {\em plus} the first
Betti number of the graph. Observe that both relations preserve this new
degree.

Read backwards, our results give an interpretation of the
Kontsevich integral in terms of the geometrically defined
equivalence relations of grope cobordism.

\section{Grope concordance} \label{sec:4d}
\setzero\vskip-5mm \hspace{5mm}

\begin{figure}[ht]\label{fig:symmetric gropes}
\begin{center}
\includegraphics[scale=.5]{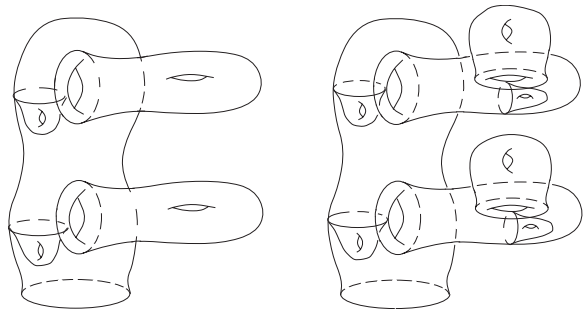}
\caption{Symmetric gropes of height~$2$ and $2.5$}
\end{center}
\end{figure}

We now turn to the 4-dimensional aspects of the theory.
It may look like the end of the story to realize that any knot with trivial
Arf invariant bounds a grope of arbitrary big class embedded in $D^4$,
\cite[Prop.3.8]{CT2}. However, group theory has more to offer than the
lower central series. Recall that the {\em derived series} of a group $G$
is defined inductively by the iterated commutators
$$
G^{(1)}:=[G,G] \text{ and } G^{(h)}:=[G^{(h-1)},G^{(h-1)}] \text{ for }
h>1.
$$
Accordingly, we may define {\em symmetric gropes} with their complexity now
measured by {\em height}, satisfying the following defining property: A
continuous map
$\phi:S^1\to X$ represents an element in $\pi_1X^{(h)}$  if
and only if it extends to a continuous map of a symmetric grope of
height~$h$. Thus a symmetric grope of height~1 is just a surface, and a
symmetric grope of height~$h$ is obtained from a bottom surface by
attaching symmetric gropes of height~$(h-1)$ to a {\em full} symplectic
basis of curves. This defines symmetric gropes of height~$h$ with one
boundary circle as in Figure~\ref{fig:symmetric gropes}

Note that a symmetric grope of height~$h$ is also a grope of class~$2^h$,
just like in group theory. But conversely, not every grope is symmetric.
It should also be clear from
Figure~\ref{fig:symmetric gropes} how one defines symmetric gropes with
half-integer height (even though there is no group theoretic analogue).

In the following definition we attempt to distinguish the terms
``cobordant'' and ``concordant'' in the sense that the latter refers to 4
dimensions, whereas the former was used in dimension 3, see
Definition~\ref{def:3d}. Historically, these terms were used
interchangeably, but we hope not to create any confusion with our new
distinction.

\begin{definition} \label{def:4d}
 Two oriented knots in $S^3$ are {\em grope concordant}
of height~$h\in\frac{1}{2}\N$, if there is an embedded symmetric grope of
height~$h$ in $S^3\times [0,1]$ such that its boundary consists exactly of
the given knots $K_i:S^1\hra S^3 \times\{i\}$.
\end{definition}

Observe that since an annulus is a symmetric grope of arbitrary height we
indeed get a filtration of the {\em knot concordance group}. This group is
defined by identifying two knots which cobound an embedded annulus in
$S^3\times [0,1]$, where there are two theories depending on whether the
embedding is smooth or just topological (and locally flat). For grope
concordance the smaller topological knot concordance group is the more
natural setting: a locally flat
topological embedding of a grope (defined by the obvious local model at
the singular points) can be perturbed to become smooth. This perturbation
might introduce many self-intersection points in the surface stages of the
grope. However, these new singularities are arbitrarily small and thus
they can be removed at the expense of increasing the genus of the surface
stage in question  but without changing the height of the grope.

In joint work with Tim Cochran and Kent Orr
\cite{COT1}, we showed that all known knot concordance invariants fit
beautifully into the scheme of grope concordance! In particular, all known
invariants turned out to already be invariants of grope concordance of
height~$3.5$:

\begin{theorem} \label{thm:4d}
Consider two knots $K_i$ in  $S^3$. Then
\begin{enumerate}
\item $K_i$ have the same Arf invariant if and only if they are grope
concordant of height~$1.5$ (or class~3).
\item $K_i$ are {\em algebraically concordant} in the sense of Levine
\cite{Le} (i.e.\ all twisted signatures and twisted Arf invariants agree)
if and only if they are grope concordant of height~$2.5$.
\item If $K_i$ are grope concordant of height~$3.5$ then they have the same
Casson-Gordon invariants \cite{CG}.
\end{enumerate}
\end{theorem}

The third statement includes the generalizations of Casson-Gordon
invariants by Gilmer, Kirk-Livingston, and Letsche.
In \cite{COT1} we prove an even stronger version of the third part of
Theorem~\ref{thm:4d}. Namely, we give a weaker condition for a knot $K$ to
have vanishing Casson-Gordon invariants: it suffices that $K$ is
{\em $(1.5)$-solvable}. All the obstruction theory in \cite{COT1} is based
on the definition of {\em $(h)$-solvable} knots, $h\in\frac{1}{2}\N$, which
we shall not give here. Suffice it to say that this definition is closer to
the algebraic topology of $4$-manifolds than grope concordance. In
\cite[Thm.8.11]{COT1} we show that a knot which bounds an embedded grope
of height $(h+2)$ in $D^4$ is $(h)$-solvable.

 It should come as no surprise that the invariants which
detect grope concordance have to do with solvable quotients of the knot
group. In fact, the above invariants are all obtained by studying the Witt
class of the intersection form of a certain $4$-manifold $M^4$ whose
boundary is obtained by $0$-surgery on the knot. The different cases are
distinguished by the fundamental group $\pi_1M$, namely
\begin{enumerate}
\item $\pi_1M$ is trivial for the Arf invariant,
\item $\pi_1M$ is infinite cyclic for algebraic concordance, and
\item $\pi_1M$ is a dihedral group for Casson-Gordon invariants.
\end{enumerate}
So the previously known concordance invariants stopped at solvable groups
which are extensions of abelian by abelian groups. To proceed further in
the understanding of grope (and knot) concordance, one must be able to
handle more complicated solvable groups.
A program for that purpose was developed in \cite{COT1} by giving an
elaborate boot strap argument to construct inductively representations of
knot groups into certain {\em universal solvable} groups. On the way,
we introduced Blanchfield duality pairings in solvable covers of
the knot complement by using noncommutative localizations of the group
rings in question.

The main idea of the boot strap is that a particular choice of
``vanishing'' of the previous invariant {\em defines} the map into the
next solvable group (and hence the next invariant). In terms of gropes this
can be expressed quite nicely as follows: pick a grope concordance of
height~$h\in\N$ and use it to construct a certain
$4$-manifold whose intersection form gives an obstruction to being able to
extend that grope to height~$h.5$. There is an obvious technical problem
in such an approach, already present in \cite{CG}: to show that there is
no grope concordance of height~$h.5$, one needs to prove
non-triviality of the obstruction {\em for all} possible gropes of
height~$h$. One way around this problem is to construct examples where the
grope concordances of small height are in some sense unique. This was done
successfully in \cite{COT1} for the level above Casson-Gordon invariants,
and in \cite{COT2} we even obtain the following infinite generation result.
Let $\cG_h$ be the graded quotient groups of knots, grope concordant of
height~$h$ to the unknot, modulo grope concordance of height~$h.5$. Then
the results of Levine and Casson-Gordon show that $\cG_2$ and $\cG_3$ are
not finitely generated.

\begin{theorem} \label{thm:cot2}
$\cG_4$ is not finitely generated.
\end{theorem}

The easiest example of a non-slice knot with vanishing Casson-Gordon
invariants is given in \cite[Fig.6.5]{COT1}.
As explained in the introduction, the last step in the proof of
Theorem~\ref{thm:cot2} is to show that the intersection form of the
$4$-manifold in question is nontrivial in a certain Witt group. Our new
tool is the von Neumann signature which has the additional bonus that it
takes values in
$\R$, which is not finitely generated as an abelian group. This fact makes
the above result tractable. We cannot review any aspect of the von Neumann
signature here, but see \cite[Sec.5]{COT1}.

Last but not least, it should be mentioned that we now know that for every
$h\in\N$ the groups $\cG_h$ are nontrivial. This work in progress
\cite{CoT3} uses as the main additional input the Cheeger-Gromov estimate
for von Neumann $\eta$-invariants \cite{ChG} in order to get around the
technical problem mentioned above. It is very likely that non of the groups
$\cG_h, h\in\N, h\geq 2,$ are finitely generated.

\label{lastpage}

\end{document}